\documentclass{article}

\usepackage{amssymb}

\newtheorem{theorem}{Theorem}[section]

\newtheorem{e-proposition}[theorem]{Proposition}

\newtheorem{corollary}[theorem]{Corollary}

\newtheorem{e-definition}[theorem]{Definition\rm}

\newtheorem{remark}{\it Remark\/}

\newtheorem{example}{\it Example\/}

\def\og{\leavevmode\raise.3ex\hbox{$\scriptscriptstyle\langle\!\langle$~}}

\def\fg{\leavevmode\raise.3ex\hbox{~$\!\scriptscriptstyle\,\rangle\!\rangle$}}

\begin{document}

\begin{center}

\vspace{-4.0cm}

{\Large \bf {\rm BMO} is the intersection of two}
\vspace{0.5cm}

{\Large \bf translates of dyadic {\rm BMO}}
\vspace {0.4cm}

{\large Tao MEI}

Dept of Math, TAMU, Tx 77840, USA, tmei@math.tamu.edu;

\begin{abstract}
Let $\mathbb T$ be the unit circle
on ${\Bbb R}^2$. Denote by {\rm BMO}$({\Bbb T})$ the classical {\rm BMO}
space
and denote by {\rm BMO}$_{\mathcal D}({\Bbb T})$ the usual dyadic BMO space
on ${\Bbb T}$. Then, for suitably chosen $\delta \in {\Bbb R},$
we have
\[
\left\| \varphi \right\| _{{\rm BMO}({\Bbb T})}\backsimeq \left\|
\varphi \right\| _{{\rm BMO}_{{\mathcal D}}({\Bbb T})}+\left\| \varphi
(\cdot -2\delta \pi )\right\| _{{\rm BMO}_{{\mathcal D}}({\Bbb T}
)},\forall \varphi \in {\rm BMO({\Bbb T})}
\]
{\it To cite this article:    C. R.
Acad. Sci. Paris, Ser. I 336 (2003).}

\vskip 0.5\baselineskip

\end{abstract}
\end{center}

\section{Introduction}

Let $\Bbb{T}$ be the unit circle on $\Bbb{R}^2$,
identified with $(0,2\pi]$.
 Recall that
\[
\mathrm{BMO}(\Bbb{T})=\{\varphi \in L^1(\Bbb{T}):\left\| \varphi
\right\|_{
\mathrm{BMO}(\Bbb{T})}=\sup \{\frac 1{|I|}\int_I|\varphi -\varphi
_I|d\theta
\}<\infty \}
\]
where the supremum runs over all intervals $I$ on $\Bbb{T}$ and $\varphi
_I=\frac 1{|I|}\int_I\varphi (s)ds.$ Let $\mathcal{D}=\{\mathcal{D}
_n\}_{n\geq 0}$ be the family of the usual dyadic $\sigma -$algebra on
$\Bbb{ T},$ i.e.
\[
\mathcal{D}_n=\sigma \{(D_n^k)_{0\leq k<2^n}\},\quad D_n^k=({2\pi
 }k{2^{-n}
}, {2\pi }(k+1){2^{-n}}];n\geq 0.
\]
Recall that the usual dyadic BMO space is defined by
\[
\mathrm{BMO}_{\mathcal{D}}(\Bbb{T})=\{\varphi \in L^1(\Bbb{T}):\left\|
\varphi \right\| _{\mathrm{BMO}_{\mathcal{D}}(\Bbb{T})}=\sup_{n,k}\{\frac
{2^n
}{2\pi }\int_{D_n^k}|\varphi -\varphi _{D_n^k}|d\theta \}<\infty \}.
\]
$\mathrm{BMO}(\Bbb{T})$ and the dyadic BMO space
 $\mathrm{BMO}_{\mathcal{D}
}(
\Bbb{T})$ have many similarities, but nevertheless certain differences.
The dyadic
\textrm{BMO} space is usually much easier to study. Some works have been
done to study the relationship between the two kinds of \textrm{BMO
}spaces (see
\cite{GJ},
\cite{p}). In this paper, we show that, for any positive $\delta $ 
suitably
chosen (more precisely satisfying $d(\delta )>0,$ with $d(\delta )$ as
defined below) $\varphi $ is in $\mathrm{BMO}(\Bbb{T})$ if and only if $
\varphi (\cdot )$ and $\varphi (\cdot -2\pi\delta )$ are in
$\mathrm{BMO}_{
\mathcal{D}}(\Bbb{T}).$ Clearly the analogous result holds 
on $\Bbb{R}$ with
the same proof (see the final remark below).

\section{The Main Result}

Let $A$ be the collection of all dyadic rationals.
For $0<\delta <1,$ define
its relative distance to $A,$ denoted by $d(\delta )$ in this paper, as
follows
\[
d(\delta ):=\inf \{2^n|\delta  - k{2^{-n}}| \mid n\geq 0,k\in \Bbb{Z}\}.
\]
Let $\Bbb{T}$ be the unit circle in $\Bbb{R}^2.$ For $\delta $ with $
d(\delta )>0,$ we consider the filtration $\mathcal{D}^\delta =\{\mathcal
{D}
_n^\delta \}_{n\geq 0}\,\,$ on $\Bbb{T}$ obtained from the usual dyadic
filtration after translation by $2\pi\delta$. More precisely:
\[
\mathcal{D}_n^\delta =\sigma \{(D_n^{\delta ,k})_{0\leq
k<2^n}\}, \ D_n^{\delta ,k}=(2\delta \pi +  {2\pi }k{2^{-n}},2\delta \pi
+{2\pi }(k+1){2^{-n}} ],\ \forall n\geq 0.
\]
Hence, if we define $\left\| \varphi \right\|
_{\mathrm{BMO}_{\mathcal{D}^\delta }(\Bbb{T}))}$ in the usual way, we
have
$
\left\| \varphi \right\| _{\mathrm{BMO}_{\mathcal{D}^\delta }(\Bbb{T}
))}=\left\| \varphi (\cdot -2\delta \pi )\right\|
_{\mathrm{BMO}_{\mathcal{D} }(\Bbb{T})}.
$

In this paper, we will say $\mathcal{D}$
 (resp. $\mathcal{D}^\delta$) ``fits"
an interval $I\subset \Bbb{T}$ with 
 fit-constant $c$ if there exist $n\geq 0,0\leq k_I<2^n$
such that $I\subset $ $D_n^{k_I}$ (resp. $I\subset $ $D_n^{\delta ,k_I}$)
 and $|D_n^{k_I}|\leq c|I|$ (resp. $|D_n^{\delta ,k_I}|\leq c|I|$).
Our key observation is the following simple fact.

\begin{e-proposition}
For any interval $I\subset \Bbb{T}$, either
 $\mathcal {D}$ or $\mathcal {D}_\delta$
fits $I$ with
fit-constant $ 2 /{d(\delta)}$.
\end{e-proposition}

\textbf{Proof} If $|I|\geq 2\pi d(\delta ),$ let $n=k_I=0,$ then $I\subset
D_0^0=(0,2\pi ].$

If $|I|<2\pi d(\delta ),$ let $n\geq 0$
 be the integer such that $d(\delta )
 {2\pi }{2^{-n-1}}\leq |I|<d(\delta ){2\pi }{2^{-n}}.$ Set
\[A_n=\{k {2\pi }{2^{-n}}; 0\leq k<2^n\},\quad
 A_n^\delta =\{2\delta \pi +k
 {2\pi }{2^{-n}};0\leq k<2^n\Bbb{\}}.
\]
Note that for any two points $a,b\in A_n\cup A_n^\delta ,$ we have $
|a-b|\geq d(\delta ) {2\pi }{2^{-n}}>|I|.$ Thus there is at most one
element of $A_n\cup A_n^\delta $ belonging to $I.\,$ Then $I\cap A_n=\phi
$ or $I\cap A_n^\delta=   \phi .$ Therefore, $I$ must be contained in
some
$ D_n^{k_I}$ or $D_n^{\delta ,k_I}$ and $|D_n^{k_I}|=|D_n^{\delta ,k_I}|=
\frac {2\pi}{2^n}
\leq 2/{d(\delta
)}|I|.$     
\begin{remark}

From the above proposition, a number of ``classical" results become
immediate consequences of their ``probabilistic" counterparts. For
instance, Doob's maximal inequality implies the Hardy-Littlewood maximal
inequality immediately. 

\end{remark}

\begin{theorem}
For $\varphi \in L^1(\Bbb{T}),0<\delta <1$ with $d(\delta )>0,$ we have
\[
\left\| \varphi \right\| _{\mathrm{BMO}(\Bbb{T})}\leq \frac 4{d(\delta
)}\max \{\left\| \varphi \right\| _{\mathrm{BMO}_{\mathcal{D}}(\Bbb{T}
)},\left\| \varphi (\cdot -2\delta \pi )\right\|
_{\mathrm{BMO}_{\mathcal{D }
}(\Bbb{T})}\}.
\]
\end{theorem}

\textbf{Proof}
By the  above proposition, for every interval $I\subset {\Bbb T}$, there
exist
$N,k_{I}$ such that $I\subseteq D_N^{k_{I}}$ or $I\subseteq D_N^{\delta
,k_I}$ and
$\frac {2\pi} {2^N}\leq \frac 2{d(\delta )}|I|.$ If $D_N^{\delta ,k_{I}}$
contains $I,$ then
\begin{eqnarray*}
\frac 1{|I |}\int_{I }|\varphi (\theta)-\varphi
_{I }|d\theta &\leq &\frac 1{|I
|}\int_{I }|\varphi (\theta)-\varphi _{D_N^{\delta ,k_{I
}}}|d\theta +|\varphi _{D_N^{\delta ,k_{I }}}-\varphi
_{I }| 
\leq \frac 2{|I |}\int_{I }|\varphi (\theta)-\varphi
_{D_N^{\delta ,k_{I }}}|d\theta\\
&\leq &\frac 4{d(\delta )|D_N^{\delta ,k_{I
}}|}\int_{D_N^{\delta ,k_{I }}}|\varphi (\theta)-\varphi
_{D_N^{\delta ,k_{I }}}|d\theta 
\leq\frac 4{d(\delta )}\left\| \varphi \right\|
_{\mathrm{BMO}_{\mathcal{D }^{^\delta }}(\Bbb{T})} \\
&&=\frac 4{d(\delta )}\left\| \varphi (\cdot -2\delta \pi )\right\| _{
\mathrm{BMO}_{\mathcal{D}}(\Bbb{T}).}
\end{eqnarray*}
If $D_N^{k_{I }}$ contains $I ,$ then similarly
\[
\frac 1{|I |}\int_{I }|\varphi (\theta)-\varphi
_{I }|d\theta \leq \frac 4{d(\delta )}\left\| \varphi \right\| _{
\mathrm{BMO}_{\mathcal{D}}(\Bbb{T})}
\]
Thus, taking the supremum over all intervals $I\subset {\Bbb T} $, we get
\[
\left\| \varphi \right\| _{\mathrm{BMO}(\Bbb{T})}\leq \frac 4{d(\delta
)}\max \{\left\| \varphi \right\| _{\mathrm{BMO}_{\mathcal{D}}(\Bbb{T}
)},\left\| \varphi (\cdot -2\delta \pi )\right\| _{\mathrm{BMO}_{\mathcal{D}}(\Bbb{T})}\}.
\ \ \ \ \ \ \ \ \ \
\]

\begin{example}
\textbf{\ }Let $\delta=   1/3,$ then $d(\delta )=  1/3,$ and then
\[
\left\| \varphi \right\| _{\mathrm{BMO}(\Bbb{T})}\leq 12\max \{\left\|
\varphi \right\| _{\mathrm{BMO}_{\mathcal{D}}(\Bbb{T})},\left\| \varphi
(\cdot -\frac{2\pi }3)\right\| _{\mathrm{BMO}_{\mathcal{D}}(\Bbb{T})}\}.
\]
\end{example}

\begin{remark}
Let $\varphi ^{\#}(t)=\sup_{I\ni t}\frac
1{|I|}\int_I|\varphi-\varphi_I|d\theta$ and
$\varphi^{\#}_\mathcal {D}(t)=\sup_{D_n^k\ni t}\frac 1{|D_n^k|}\int_{D_n^
k}|\varphi-\varphi
_{D_n^k}|d\theta$.
It is easy to see that $\{\delta, d(\delta)>0\}$ is exactly
the set of all $\delta$'s such that $\varphi^{\#}\leq
c\max\{\varphi^{\#}_\mathcal {D} ,\varphi^{\#}_\mathcal
{D}(\cdot-2\pi\delta)\}$ for some $c>0$. The same statement trivially
remains valid in the Banach space valued case and is particularly useful
in the operator valued case: see \cite{m} for some results in that
direction.
\end{remark}

\begin{remark}
One can check that the set $\{\delta, d(\delta)>0\}$
is dense in $(0,1)$ while its measure is zero.
\end{remark}

\begin{corollary}
$\mathrm{BMO} (\Bbb{T})=\mathrm{BMO}_{\mathcal{D}}(\Bbb{T})\cap
\mathrm{BMO}_{
\mathcal{D^\delta }}(\Bbb{T})$ with equivalent norms.
\end{corollary}

Denote by $H_{\mathcal{D}}^1$ (resp. $H_{\mathcal{D}^\delta }^1$)
the dyadic
Hardy space with respect to $\mathcal{D}$ (resp. $\mathcal{D}^\delta $).
 By
duality, we have

\begin{corollary}
$H^1=H_{\mathcal{D}}^1+H_{\mathcal{D}^\delta }^1$ with equivalent
norms.
\end{corollary}

\begin{remark}
There is another way to see Corollary 2.4.
Denote by $H^{1,at}$ the classical atomic Hardy space. Denote by
$H^{1,at}_{\mathcal {D}}$ (resp. $H_{\mathcal{D}^\delta }^{1,at}$) the
dyadic atomic Hardy space with respect to $\mathcal{D}$ (resp.
$\mathcal{D}^\delta $). From Proposition 2.1, we see that any atom is a
dyadic atom (up to a fixed factor) with respect to either $\mathcal{D}$
or $\mathcal{D}^\delta$. Thus
$H^{1,at}=H_{\mathcal{D}}^{1,at}+H_{\mathcal{D}^\delta }^{1,at}$ with
equivalent norms. Since $H^{1,at}=H^1$ and $H^{1,at}_{\mathcal
{D}}=H_{\mathcal{D}}^1$, we obtain Corollary 2.4.
\end{remark}

\begin{remark}
See \cite{p} for a recent result (of the same flavor) comparing Hilbert
transforms and martingale transforms proved by averaging 
shifted and dilated
dyadic filtrations.
\end{remark}

\begin{remark}
John Garnett kindly informed us that he already knew that ${\rm BMO}
(\Bbb{T})$ coincides with the intersection of three
(suitably chosen) translates of dyadic
{\rm BMO}$(\Bbb{T})$
 (the
idea for this can be traced back to page 417 of \cite{G}), but our main
result seems new.
\end{remark}

We now turn to  the case of   dimension $m>1$.
By a straightforward product argument, one can deduce
from the above proposition that $\mathrm{BMO}(\Bbb{T}^m)$
coincides with the intersection of a family
of $2^m$ translates of the dyadic version of $\mathrm{BMO}(\Bbb{T}^m)$.
However, we wish to show below
that the number of translates can be reduced to $m+1$.

In the following, we always suppose $\{\delta_i\}_{i=0}^m$ is a sequence
in $(0,1)$ such that $$d(\{\delta_i\}_{i=0}^m):=\min_{i\neq
j}d(\delta_i-\delta_j)> 0.$$
Let $\mathcal {D}^{\delta_i}$ be the translation by
 $2\pi\delta_i$ of the family of the usual (one dimensional) dyadic
$\sigma$-algebra. Set $\mathcal {F}^i=(\mathcal {D}^{\delta_i})^m, 0\leq
i\leq m$. Then we  get $m+1$
families of increasing dyadic $\sigma $-algebras on $\Bbb T^m$.

\begin{e-proposition}
Let $\mathcal {F}^i, 0\leq i\leq m$ be as above and let
$c=2/d(\{\delta_i\}_{i=0}^m)$. Then, for any cube $J\subset {\Bbb T}^m$,
there exists some $\mathcal {F}^i$
which fits $J$ with  fit-constant   $c^m.$
\end{e-proposition}

{\bf Proof} Write $J\subset \Bbb T^m$ as $J=J_1\times J_2 \times
\cdot\cdot\cdot \times J_m$, where
$J_i$ are intervals in ${\Bbb T}, 1\leq i\leq m$. Let
$\{\delta_i\}_{i=0}^m$ be such that $d(\{\delta_i\}_{i=0}^m)>0$.
By Proposition 2.1, for every $J_i$,
there is at most one $k_i, 0\leq k_i \leq m$ such that $\mathcal
{D}^{\delta_{k_i}}$
 does not fit $J_i$ with constant $c$.
Then there is at least one $\mathcal {D}^{\delta_k}$ which 
fits all $J_i$ with constant $c$.
Thus (with an obvious extension of our terminology)
we may say that $\mathcal {F}^k$ fits $J$ with fit-constant $c^m$. 
From Proposition 2.5 we have
\begin{theorem}
(In the case of $\Bbb{T}^m$)
Let $\{\delta_i\}_{i=0}^m$ be
a sequence in $(0,1)$ such that $d(\{\delta_i\}_{i=0}^m)>0$.
Let $^i\delta=(\delta_i, \delta_i,\cdot\cdot,\delta_i)$.
Then, for $\varphi\in L^1(\Bbb{T}^m)$, we have
\[
\left\| \varphi \right\| _{\mathrm{BMO}(\Bbb{T}^m)}\leq 2(
2/{d(\{\delta_i\}_{i=0}^m)})^m
\max_{0\leq i\leq m}\{\left\| \varphi (\cdot -{^i\delta}2\pi)\right\|_
{\mathrm{BMO}_{\mathcal{D}}(\Bbb{T}^m)}\}.
\]
\end{theorem}

\begin{remark}
To extend our results to $\Bbb{R}^m,$ denote
 by $\mathcal{D}(\Bbb{R})$ the
family of the usual dyadic $\sigma-$algebra on $\Bbb{R}.$
For $0<\delta <1$
with $d(\delta )>0,$ choose an increasing family of
 dyadic $\sigma-$algebra
$\mathcal{D}^\delta (\Bbb{R})=(\mathcal{D}_n^\delta )_{n\in
\Bbb{Z}}(\Bbb{R})
$ such that, for $n$ even,
\[
\mathcal{D}_n^\delta (\Bbb{R})=\sigma (\{D_n^{\delta ,k}\}_{k\in
\Bbb{Z}}),
\begin{array}{l}
D_n^{\delta ,k}(\Bbb{R})=(\frac k{2^n}+\delta ,\frac{k+1}{2^n}+\delta
], \ n\geq 0, \\
D_n^{\delta ,k}(\Bbb{R})=(\frac k{2^n}+\delta +\sum_{j=n+2}^0\frac
1{2^j},\frac{k+1}{2^n}+\delta +\sum_{j=n+2}^0\frac 1{2^j}], \
n<0.
\end{array}
\]
Note that all $\mathcal{D}^{\delta}_n(\Bbb{R})$'s are given after fixing $\mathcal{D}_n^{\delta}(\Bbb{R})$'s
for all even $n$'s.
Let $\{\delta_i\}_{i=0}^m$ be
a sequence in $(0,1)$ such that $d(\{\delta_i\}_{i=0}^m)>0$.
Let $^i\mathcal{D}^\delta (\Bbb{R}^m)=\{^i\mathcal{D}^\delta_n(\Bbb{R}^m)
\}_{n\in \Bbb {N}}$, where
$^i\mathcal{D}^\delta_n(\Bbb{R}^m)$ is the $m$ times product of
the $\sigma$-algebra
$\mathcal{D}^{\delta_i}_n(\Bbb{R})$.
Then,
by the same idea as above, we can get
\[
\left\| \varphi \right\| _{\mathrm{BMO}(\Bbb{R}^m)}
\leq 2(4/{d(\{\delta_i\}_{i=0}^m)})^m
\max_{0\leq i\leq m}(\left\| \varphi
\right\|_{\mathrm{BMO}_{^i\mathcal{D}^\delta }(\Bbb{R}^m)})\quad \forall
\varphi\in  L^1(\Bbb{R}^m) .
\]
\end{remark}

\section*{Acknowledgements}
The author is very grateful to Q. Xu and his adviser G. Pisier for useful
conversations.

\end{document}